\documentstyle[amssymb,12pt]{article}
\def\Tr{\mathop{\,\mathrm{Tr}\,}\nolimits}
\def\trace{\mathop{\,\mathrm{tr}\,}\nolimits}
\def\nn#1{{\bf \#{#1}}}
\def\rom#1{{\normalshape #1}}
\newenvironment{mat}{\footnotesize$\left[\begin{array}}
                    {\end{array}\right]$}
\newtheorem{theorem}{Theorem}
\newenvironment{proof}{\begin{trivlist}\item[]
{\bf Proof.} }{\hfill $\Box$ \end{trivlist}}
\begin{document}
\begin{center}
{\large\bf HOMOGENEOUS HYPERSURFACES WITH}\\[8pt]
{\large\bf ISOTROPY IN AFFINE FOUR-SPACE}\\[20pt]
{\large\sc Michael Eastwood and Vladimir Ezhov}\end{center}
\def\thefootnote{}
\footnotetext{The authors gratefully acknowledge support from the Australian
Research Council.}
\def\thefootnote{\fnsymbol{footnote}}\setcounter{footnote}{1}
{\small {\bf Abstract:} We classify the non-degenerate homogeneous
hypersurfaces in real and complex affine four-space whose symmetry group
is at least four-dimensional.}

\section{Statement of results}\label{results}
\renewcommand{\arraystretch}{1.2}
\begin{center}\begin{tabular}{|l|l|l|l|c|}\cline{2-5}
\multicolumn{1}{c|}{}&Equation & Basepoint & Parameter    & Dimension \\
\multicolumn{1}{c|}{}&         &           & Restrictions & of Isotropy\\
\hline\hline
\nn{1} &$W=XY+Z^2$                &$(0,0,0,0)$&                     &4\\ \hline
\nn{2} &$W^2=XY+Z^2+1$            &$(1,0,0,0)$&                     &3\\ \hline
\nn{3} &$W=XY+Z^2+X^3$            &$(0,0,0,0)$&                     &2\\ \hline
\nn{4} &$W=XY+Z^2+X^2Z+\alpha X^4$&$(0,0,0,0)$&$\alpha$ arbitrary   &1\\ \hline
\nn{5} &$W=XY+Z^2+XZ^2$           &$(0,0,0,0)$&                     &1\\ \hline
\nn{6} &$W^2=XY+X^2Y+X^2Z$        &$(1,1,0,1)$&                     &1\\ \hline
\nn{7} &$W=XY+Z^\alpha$           &$(1,0,0,1)$&$\alpha\not=0,1,2$   &1\\ \hline
\nn{8} &$W=XY+e^Z$                &$(1,0,0,0)$&                     &1\\ \hline
\nn{9} &$W=XY+\log Z$             &$(0,0,0,1)$&                     &1\\ \hline
\nn{10}&$W=XY+Z\log Z$            &$(0,0,0,1)$&                     &1\\ \hline
\nn{11}&$W^2=XY+Z^\alpha$         &$(1,0,0,1)$&$\alpha\not=0,1,2$   &1\\ \hline
\nn{12}&$W^2=XY+e^Z$              &$(1,0,0,0)$&                     &1\\ \hline
\nn{13}&$WZ=XY+Z^\alpha$          &$(1,0,0,1)$&$\alpha\not=0,1,2$   &1\\ \hline
\nn{14}&$WZ=XY+Z\log Z$           &$(0,0,0,1)$&                     &1\\ \hline
\nn{15}&$WZ=XY+Z^2\log Z$         &$(0,0,0,1)$&                     &1\\ \hline
\nn{16}&$W=XY+Z^2+X^\alpha$       &$(1,1,0,0)$&$\alpha\not=0,1,2,3$ &1\\ \hline
\nn{17}&$W=XY+Z^2+e^X$            &$(1,0,0,0)$&                     &1\\ \hline
\nn{18}&$W=XY+Z^2+\log X$         &$(0,1,0,0)$&                     &1\\ \hline
\nn{19}&$W=XY+Z^2+X\log X$        &$(0,1,0,0)$&                     &1\\ \hline
\nn{20}&$W=XY+Z^2+X^2\log X$      &$(0,1,0,0)$&                     &1\\ \hline
\end{tabular}\end{center}\renewcommand{\arraystretch}{1}
Each of the equations in this table defines, near its basepoint, a
non-degenerate homogeneous hypersurface in complex affine four-space. This
means that it may be analytically continued to an orbit of a Lie subgroup of
the group of affine symmetries. In each case, the full symmetry group (the
maximal such Lie subgroup) has dimension at least four. In other words, there
is a non-trivial (positive dimension) Lie subgroup preserving the basepoint.
Indeed, the dimension of this {\em isotropy\/} is as listed. Furthermore, this
is a complete list:
\begin{theorem}\label{explicit} Every homogeneous hypersurface with isotropy in
complex affine four-space may, for a suitable choice of affine co\"ordinate
system, be found in the table above. Different entries in this table and
different values of the parameter $\alpha$ define affinely distinct
hypersurfaces.\end{theorem}

The corresponding real list is given at the end of this article.
The complex classification list for non-degenerate surfaces in affine
three-space has just three entries, namely the graph of a non-degenerate
quadratic, the complex sphere, and the Cayley surface. This is proved
explicitly in~\cite{dkr}. Theorem~\ref{explicit} will be deduced from an
alternative formulation as follows.
\begin{theorem}\label{maintheorem}
Every homogeneous hypersurface with isotropy in complex affine four-space may
be found in the following list of normal forms for a suitable choice of affine
co\"ordinates and free parameter~$b$:
\begin{description}\addtolength{\itemsep}{-7pt}
                   \setlength{\labelwidth}{30pt}
                   \setlength{\itemindent}{0pt}
\item[\rom{Qd}\hfill]   $w=2xy+z^2+{\mathrm{O}}(5)=2xy+z^2$
\item[\rom{Sp}\hfill]   $w=2xy+z^2+4x^2y^2+4xyz^2+z^4+{\mathrm{O}}(5)
                          =(1-\sqrt{1-4(2xy+z^2)})/2$
\item[\rom{I3}\hfill]   $w=2xy+z^2+x^3+{\mathrm{O}}(5)=2xy+z^2+x^3$
\item[\rom{I2}\hfill]   $w=2xy+z^2+x^2z+bx^4+{\mathrm{O}}(5)=2xy+z^2+x^2z+bx^4$
\item[\rom{I1.1}\hfill] $w=2xy+z^2+x^2y-2xz^2+\frac{1}{2}x^3y
                         -x^2z^2+{\mathrm{O}}(5)
                          =\displaystyle\frac{4xy+2z^2-5xz^2}{2-x}$
\item[\rom{I1.2}\hfill] $w=2xy+z^2+x^2y-2xz^2+\frac{1}{2}x^3y
                         +\frac{21}{4}x^2z^2+{\mathrm{O}}(5)
                         =\displaystyle\frac{4(2xy+z^2+5x^2y)}{(2-x)(2+5x)}$
\item[\rom{I0.1}\hfill] $w=2xy+z^2+3xyz-z^3
     +\frac{9}{4}x^2y^2-\frac{9}{2}xyz^2+\frac{15}{32}bz^4+{\mathrm{O}}(5)$
\item[\rom{I0.2}\hfill] $w=2xy+z^2+3xyz-z^3
     -\frac{9}{8}(5b-7)x^2y^2+\frac{9}{8}(5b-9)xyz^2+\frac{15}{32}bz^4
     +{\mathrm{O}}(5)$
\item[\rom{I0.3}\hfill] $w=2xy+z^2+3xyz-z^3$\hfill\mbox{ }\linebreak\mbox{ }
     $-\frac{1}{4}(b+1)(b+7)x^2y^2
      -\frac{1}{4}(b^2-7b-26)xyz^2-\frac{1}{16}(b^2-7b-14)z^4
+{\mathrm{O}}(5)$
\item[\rom{Inr}\hfill]  $w=2xy+z^2+x^3+x^4+bx^5+{\mathrm{O}}(6)$
\end{description}
In each case, the higher order terms are determined by the specified terms and
no lower order truncation has this property. These hypersurfaces are affinely
distinct save for the following three overlaps:
\begin{itemize}\addtolength{\itemsep}{-7pt}
\item Case\/ {\bf \rom{I0.1}} and case\/ {\bf \rom{I0.2}} agree when $b=1$
\item Case\/ {\bf \rom{I0.1}} and case\/ {\bf \rom{I0.3}} agree when $b=-4$
\item Case\/ {\bf \rom{I0.2}} and case\/ {\bf \rom{I0.3}} agree when $b=7/2$.
\end{itemize}
\end{theorem}
This list is obtained by choosing co\"ordinates so that the Taylor series of
the defining function and, consequently, the isotropy of the surface is in a
preferred form. In principle, this gives an algorithm for locating the surface
and its parameter as a local invariant. In cases~{\bf I0.1}, {\bf I0.2},
and~{\bf I0.3} there is clearly some choice for what to take as parameter. For
example, the coefficient of $z^4$ might be a more natural choice in case~{\bf
I0.1}. The particular choices made are so that $b$ is continuous across the
three overlaps between these cases whilst at the same time not making any of
the formulae too unwieldy.

The following table compares Theorems~\ref{explicit} and~\ref{maintheorem}.
\renewcommand{\arraystretch}{1.2}
\begin{center}\begin{tabular}{|l|l|l|l|}\hline
Explicit & Normal & How the Parameters\\
Form & Form      &$\alpha$ and $b$ are Related\\ \hline \hline
\nn{1}  & {\bf Qd}                   &                            \\ \hline
\nn{2}  & {\bf Sp}                   &                            \\ \hline
\nn{3}  & {\bf I3}                   &                            \\ \hline
\nn{4}  & {\bf I2}                   & $\alpha=b$                 \\ \hline
\nn{5}  & {\bf I1.1}                 &                            \\ \hline
\nn{6}  & {\bf I1.2}                 &                            \\ \hline
\nn{7}  & {\bf I0.1}                 & $\alpha=(2b-2)/(b+4)$      \\ \hline
\nn{8}  & {\bf I0.1} or {\bf I0.3}   & $b=-4$                     \\ \hline
\nn{9}  & {\bf I0.1} or {\bf I0.2}   & $b=1$                      \\ \hline
\nn{10} & {\bf I0.1}                 & $b=6$                      \\ \hline
\nn{11} & {\bf I0.2}                 & $\alpha=(4b-4)/(4b-9)$     \\ \hline
\nn{12} & {\bf I0.2}                 & $b=9/4$                    \\ \hline
\nn{13} & {\bf I0.3}                 & $\alpha=15/(b+4)$          \\ \hline
\nn{14} & {\bf I0.3}                 & $b=11$                     \\ \hline
\nn{15} & {\bf I0.2} or {\bf I0.3}   & $b=7/2$                    \\ \hline
\nn{16} & {\bf Inr}                  & $\alpha=(15b-16)/(5b-4)$   \\ \hline
\nn{17} & {\bf Inr}                  & $b=4/5$                    \\ \hline
\nn{18} & {\bf Inr}                  & $b=16/15$                  \\ \hline
\nn{19} & {\bf Inr}                  & $b=6/5$                    \\ \hline
\nn{20} & {\bf Inr}                  & $b=8/5$                    \\ \hline
\end{tabular}\end{center}\renewcommand{\arraystretch}{1}

This article is organised as follows. In the next section we shall describe how
to normalise up to third order the defining function of a non-degenerate
hypersurface under the assumption that it is homogeneous with isotropy. Then,
in \S\ref{proofs} we shall use these normalisations to effect the
classification. The method follows~\cite{ee} and especially the criteria for
homogeneity developed therein. The conversion of this classification to the
list of explicit defining functions is described in~\S\ref{compare}. In
most of this
article we are working for simplicity over the complex numbers. A real
classification may be performed similarly. In \S\ref{remarks}
we make a few remarks on this task and list the real defining functions.
Though the details are different, in~\cite{lob} Loboda uses affine normal forms
to classify homogeneous surfaces in affine three-space. Presumably, his
approach could also be employed for hypersurfaces with isotropy. After our
article was completed we learned of a manuscript by N.~Mozhey who considers the
same problem with a different method.

\section{Normal forms}\label{normalforms}
We shall choose affine co\"ordinates so that the hypersurface $\Sigma$ and its
isotropy are in some preferred normal form. Firstly, we shall
choose co\"ordinates $(x,y,z,w)$ so that $\Sigma$ passes through the origin and
so that $\{w=0\}$ is its tangent plane. Recall that $\Sigma$ is supposed
non-degenerate. This allows us to normalise the quadratic terms of its defining
function. Also, we shall take the $w$-axis to be the affine normal.
The effect of these choices is that $\Sigma$ may be defined by a power series
\begin{equation}\label{series}w=F(x,y,z)=2xy+z^2+\mbox{cubic terms}
+{\mathrm{O}}(4)\end{equation}
whose cubic terms are trace-free with respect to the quadratic form
associated with $2xy+z^2$ (as explained in \cite{l}
or \cite[Proposition~1]{ee}). Specifically, this means that the cubic terms are
spanned by
\begin{equation}\label{basis}
x^3,\;x^2z,\;x^2y-2xz^2,\;3xyz-z^3,\;xy^2-2yz^2,\;y^2z,\;y^3.\end{equation}
At this stage, the remaining co\"ordinate freedom
is ${\mathrm O}(3,{\Bbb C})$ acting on $(x,y,z)$ together with the rescaling
\begin{equation}\label{rescaling}x\mapsto\lambda x\quad y\mapsto\lambda y\quad
z\mapsto\lambda z\quad w\mapsto\lambda^2w.\end{equation}
The corresponding Lie algebra may be represented by matrices of the form
\begin{equation}\label{Lie}
\mbox{\begin{mat}{cccc}t-r&0&p&0\\
                          0&t+r&-q&0\\
                            q&-p&t&0\\
                            0&0&0&2t\end{mat}}
\end{equation}
with the usual Lie bracket of matrices. As far as ${\mathrm O}(3,{\Bbb C})$ is
concerned, the adjoint action for ${\frak o}(3,{\Bbb C})$
$$\left[\rule[-22pt]{0pt}{44pt}
\mbox{\begin{mat}{ccc}-c&0&a\\ 0&c&-b\\ b&-a&0\end{mat}},
\mbox{\begin{mat}{ccc}-r&0&p\\ 0&r&-q\\ q&-p&0\end{mat}}\right]=
\mbox{\begin{mat}{ccc}aq-bp&0&ar-cp\\ 0&bp-aq&br-cq\\ cq-br&cp-ar&0
\end{mat}}$$
may be viewed as matrix multiplication
$$\mbox{\begin{mat}cp\\ q\\ r\end{mat}}\longmapsto
\mbox{\begin{mat}{ccc}-c&0&a\\ 0&c&-b\\ b&-a&0\end{mat}}
\mbox{\begin{mat}cp\\ q\\ r\end{mat}}$$
and similarly for the Adjoint action
$$\mbox{\begin{mat}cp\\ q\\ r\end{mat}}\longmapsto
M\mbox{\begin{mat}cp\\ q\\ r\end{mat}}
\quad\mbox{for }M\in{\mathrm O}(3,{\Bbb C}).$$
This standard representation has, up to scale, three orbits, namely the origin,
the vectors of non-zero length, and the non-zero null vectors. Accordingly, we
may conjugate and rescale any matrix in ${\frak o}(3,{\Bbb C})$ into one of
three standard forms:
\begin{equation}\label{threestandardforms}
\mbox{\begin{mat}{ccc}0&0&0\\ 0&0&0\\ 0&0&0\end{mat}}
\quad\mbox{or}\quad
\mbox{\begin{mat}{ccc}-1&0&0\\ 0&1&0\\ 0&0&0\end{mat}}
\quad\mbox{or}\quad
\mbox{\begin{mat}{ccc}0&0&0\\ 0&0&-1\\ 1&0&0\end{mat}}.\end{equation}
We are now in a position to normalise the isotropy of $\Sigma$. Any 1-parameter
subgroup of this isotropy will be generated by a matrix of the form~(\ref{Lie})
and the ${\frak o}(3,{\Bbb C})$ component thereof may then be normalised as
above. In the second two cases this fixes the scaling. Therefore, we obtain
three possibilities
\begin{equation}\label{generators}
\mbox{\begin{mat}{cccc}
      1&0&0&0\\ 0&1&0&0\\ 0&0&1&0\\ 0&0&0&2\end{mat}}\quad\mbox{or}\quad
\mbox{\begin{mat}{cccc}
      t-1&0&0&0\\ 0&t+1&0&0\\ 0&0&t&0\\ 0&0&0&2t\end{mat}}\quad\mbox{or}\quad
\mbox{\begin{mat}{cccc}
      t&0&0&0\\ 0&t&-1&0\\ 1&0&t&0\\ 0&0&0&2t\end{mat}}
\end{equation}
which we shall take as normal forms for an isotropy generator. The first
possi\-bi\-li\-ty occurs but is very restrictive. It means that the whole power
series~(\ref{series}) is preserved by the rescaling~(\ref{rescaling}). This
forces all cubic and higher terms to vanish. We are left with a quadratic
defining function. This is case~{\bf Qd} of Theorem~\ref{maintheorem}. Its
isotropy is four-dimensional, generated by ${\mathrm O}(3,{\Bbb C})$
and~(\ref{rescaling}). One more case which can be dealt with separately is
when all cubic terms vanish but there are some non-vanishing
higher order terms. That it must be the complex hypersphere (case~{\bf Sp}) is
an immediate consequence of the classical Maschke-Pick-Berwald Theorem (see,
e.g.~\cite{ns}) which states that a non-degenerate hypersurface with vanishing
cubic form is a hyperquadric. This classical theorem does not assume
{\em a~priori\/} that the hypersurface is homogeneous. When homogeneity is
assumed the conclusion is also an easy consequence of our approach
(see~\S\ref{proofs}).

We may now suppose that there are non-zero cubic terms in (\ref{series}) and
investigate the consequences of the corresponding hypersurface $\Sigma$
admitting either of the second two of (\ref{generators}) as a symmetry.
Writing $c(x,y,z)$ for the cubic terms, this means that
$$\left[(t-1)x\frac{\partial}{\partial x}
       +(t+1)y\frac{\partial}{\partial y}
           +tz\frac{\partial}{\partial z}\right]c(x,y,z)-2tc(x,y,z)=0$$
or
$$\left[tx\frac{\partial}{\partial x}
   +(ty-z)\frac{\partial}{\partial y}
   +(x+tz)\frac{\partial}{\partial z}\right]c(x,y,z)-2tc(x,y,z)=0,$$
respectively. Using (\ref{basis}) as a basis of the cubic terms, these
equations place $c(x,y,z)$ in the kernel of the matrices
$$\mbox{\begin{mat}{ccccccc}
t-3& 0 & 0 & 0 & 0 & 0 & 0 \\
 0 &t-2& 0 & 0 & 0 & 0 & 0 \\
 0 & 0 &t-1& 0 & 0 & 0 & 0 \\
 0 & 0 & 0 & t & 0 & 0 & 0 \\
 0 & 0 & 0 & 0 &t+1& 0 & 0 \\
 0 & 0 & 0 & 0 & 0 &t+2& 0 \\
 0 & 0 & 0 & 0 & 0 & 0 &t+3\end{mat}}$$
or
$$\mbox{\begin{mat}{ccccccc}
 t & 0 & 0 & 0 & 0 & 0 & 0 \\
 1 & t & 0 & 0 & 0 & 0 & 0 \\
 0 & -5& t & 0 & 0 & 0 & 0 \\
 0 & 0 & 3 & t & 0 & 0 & 0 \\
 0 & 0 & 0 & -2& t & 0 & 0 \\
 0 & 0 & 0 & 0 & 1 & t & 0 \\
 0 & 0 & 0 & 0 & 0 & -3& t \end{mat}},$$
respectively. The first matrix is singular if and only if $t$ is an integer in
the range $-3,\ldots,3$ whilst the second is singular only for~$t=0$. We
conclude that normalising the isotropy as we have done automatically forces
$c(x,y,z)$ to be a simple multiple of one of the seven basic
cubics~(\ref{basis}). Swopping $x$ and $y$ if necessary, we have almost proved
the following.
\begin{theorem}\label{normforms} A homogeneous hypersurface with isotropy in
complex affine four-space may be locally defined for a suitable choice of
affine co\"ordinate system by a power series of the form \rom{(\ref{series})}
which, if the cubic terms are non-zero, may be further normalised to have one
of the following forms:
\begin{description}
\item[\rom{I3}] $w=2xy+z^2+x^3+{\mathrm O}(4)$
\item[\rom{I2}] $w=2xy+z^2+x^2z+{\mathrm O}(4)$
\item[\rom{I1}] $w=2xy+z^2+x^2y-2xz^2+{\mathrm O}(4)$
\item[\rom{I0}] $w=2xy+z^2+3xyz-z^3+{\mathrm O}(4).$
\end{description}
The residual co\"ordinate freedom is generated by
\begin{equation}\label{scaling}x\mapsto\lambda^{t-1}x\quad
y\mapsto\lambda^{t+1}y\quad
z\mapsto\lambda^tz\quad w\mapsto\lambda^{2t}w\end{equation}
for $t=3,2,1,0$ respectively and, in addition,
$$\mbox{\begin{mat}c x\\ y\\ z\\ w\end{mat}}\longmapsto
\mbox{\begin{mat}{cccc} 1 & 0 & 0 & 0 \\
                      -t^2/2& 1 & -t& 0 \\
                          t & 0 & 1 & 0 \\
                          0 & 0 & 0 & 1 \end{mat}}
\mbox{\begin{mat}c x\\ y\\ z\\ w\end{mat}}$$
in case~{\bf \rom{I3}} and swopping $x$ and $y$ in case~{\bf \rom{I0}}.
\end{theorem}
\begin{proof} It is easy to verify that these normal forms do, indeed, have the
residual co\"ordinate freedoms as stated and it remains to show that this is
full extent thereof.

In case~{\bf I3} we can use the scaling freedom~(\ref{scaling}) to suppose that
$w$ is preserved on the nose, not merely up to scale. {From} the cubic term,
$x$ is then preserved up to scaling by a cube root of unity. This too may be
incorporated into (\ref{scaling}) and we may now suppose that $x$ is also
preserved on the nose. We are left with ${\mathrm O}(3,{\Bbb C})$
transformations fixing a null vector and it is easy to check that they are of
the form
\begin{equation}\label{nr}\mbox{\begin{mat}cx\\ y\\ z\end{mat}}\longmapsto
\mbox{\begin{mat}{ccc}1&0&0\\ -t^2/2&1&\mp t\\ t&0&\pm 1\end{mat}}
\mbox{\begin{mat}cx\\ y\\ z\end{mat}}.\end{equation}
With the positive sign this is a {\em null rotation\/}~\cite[p.~28]{OT}. The
negative sign may be absorbed into a scaling (\ref{scaling}) with $\lambda=-1$.
Assembling these possibilities yields precisely the freedom as stated.

In all cases $w$ and its axis are preserved up to scale and in case~{\bf I2}
the cubic term ensures that $x$ and $z$ are also preserved up to scale:
$$x\mapsto\lambda x\quad z\mapsto\nu z\quad
  w\mapsto\kappa w\quad\mbox{where }\lambda^2\nu=\kappa.$$
Now the quadratic terms force
$$y\mapsto\mu y\quad\mbox{where }\lambda\mu=\nu^2=\kappa.$$
This is of the form~(\ref{scaling}).

Writing the cubic terms in case~{\bf I1} as $x(xy-2z^2)$, a product of
irreducibles, it is clear that $x$ must be preserved up to scale as must
the quadratic form $xy-2z^2$. With the quadratic form $2xy+z^2$ also being
preserved, this easily implies that $y$ and $z$ are now preserved up to scale.
The result is of the form~(\ref{scaling}).

In case~{\bf I0}, the cubic terms factorise as $(3xy-z^2)z$ and similar
reasoning implies that $z$ and the quadratic form $xy$ are preserved,
firstly up to scale and then, by comparing scales, absolutely. The only
remaining freedom is ${\mathrm O}(2,{\Bbb C})$ acting in the $(x,y)$-variables.
The identity connected component has the form (\ref{scaling}) and the rest is
generated by the reflection which swops $x$ and~$y$.
\end{proof}

\section{Proof of Theorem~\protect\ref{maintheorem}}\label{proofs}
The proof is based on the criteria for homogeneity developed in~\cite{ee}.
For any formal power series or polynomial $G(x,y,z)$, we shall write
$\Tr^NG(x,y,z)$ for the polynomial obtained by truncation at order~$N$:
$$\mbox{if}\quad G(x,y)=\sum_{i,j,k=0}^\infty c_{ijk}x^iy^jz^k\quad
\mbox{then}\quad\Tr^NG(x,y)=\sum_{i+j+k\leq N}c_{ijk}x^iy^jz^k.$$
\begin{theorem}\label{criteria}
Suppose $f(x,y,z)$ is a polynomial of degree $N$ without constant or linear
terms. If $f(x,y,z)$ can be completed to a formal power series whose graph near
the origin is an open subset of a homogeneous hypersurface~$\Sigma$, then there
are $4\times 4$ matrices $P,Q,R$ such that
\begin{equation}\label{PQR}\raisebox{30pt}{\makebox[0pt]{$\begin{array}{rcl}
\Tr^{N-1}\mbox{\begin{mat}c
         \displaystyle\frac{\partial f}{\partial x}(x,y,z),
                      \frac{\partial f}{\partial y}(x,y,z),
                      \frac{\partial f}{\partial z}(x,y,z),-1\end{mat}}
    P\mbox{\begin{mat}cx\\ y\\ z\\ f(x,y,z)\end{mat}}
                  &=&\displaystyle-\frac{\partial f}{\partial x}(x,y,z)\\[25pt]
\Tr^{N-1}\mbox{\begin{mat}c
         \displaystyle\frac{\partial f}{\partial x}(x,y,z),
                      \frac{\partial f}{\partial y}(x,y,z),
                      \frac{\partial f}{\partial z}(x,y,z),-1\end{mat}}
    Q\mbox{\begin{mat}cx\\ y\\ z\\ f(x,y,z)\end{mat}}
                  &=&\displaystyle-\frac{\partial f}{\partial y}(x,y,z)\\[25pt]
\Tr^{N-1}\mbox{\begin{mat}c
         \displaystyle\frac{\partial f}{\partial x}(x,y,z),
                      \frac{\partial f}{\partial y}(x,y,z),
                      \frac{\partial f}{\partial z}(x,y,z),-1\end{mat}}
    R\mbox{\begin{mat}cx\\ y\\ z\\ f(x,y,z)\end{mat}}
                  &=&\displaystyle-\frac{\partial f}{\partial z}(x,y,z).
\end{array}$}}\end{equation}
Conversely, suppose that these equations have solutions $P,Q,R$ and that, for
the general such solutions,
\begin{equation}\label{closure}\Tr^{N}\mbox{\begin{mat}c
         \displaystyle\frac{\partial f}{\partial x}(x,y,z),
                      \frac{\partial f}{\partial y}(x,y,z),
                      \frac{\partial f}{\partial z}(x,y,z),-1\end{mat}}
    X\mbox{\begin{mat}cx\\ y\\ z\\ f(x,y,z)\end{mat}}=0\end{equation}
for all $X$ of the following three forms \rom{(}where $P=(p_{i,j})$
etcetera\rom{)}:
$$\begin{array}l
X=PQ-QP-(p_{1,2}-q_{1,1})P-(p_{2,2}-q_{2,1})Q-(p_{3,2}-q_{3,1})R\\[3pt]
X=QR-RQ-(q_{2,3}-r_{2,2})Q-(q_{3,3}-r_{3,2})R-(q_{1,3}-r_{1,2})P\\[3pt]
X=RP-PR-(r_{3,1}-p_{3,3})R-(r_{1,1}-p_{1,3})P-(r_{2,1}-p_{2,3})Q.
\end{array}$$
Then $f(x,y,z)$ can be uniquely completed to a formal power series whose graph
near the origin is an open subset of a homogeneous hypersurface. Furthermore,
all homogeneous hypersurfaces in affine four-space arise in this way.
\end{theorem}
\begin{proof} The proof is a simple modification of the corresponding result
for surfaces proved in Theorem~1 and Corollary~1 of~\cite{ee}. Suffice it to
say that~(\ref{PQR}), for sufficiently large~$N$, defines the symmetry algebra
of~$\Sigma$. That there are solutions is to say that there are infinitesimal
symmetries in each of the three basic co\"ordinate directions. This must be the
case if $\Sigma$ is homogeneous. For (\ref{closure}) to hold for $X$'s made out
of the general $P,Q,R$ is to say that this linear subspace of the Lie algebra
of affine motions is closed under Lie bracket. Once the symmetry algebra has
closed in this way, the higher order terms in the power series expansion of the
defining function are completely pinned down (either by exponentiating to a Lie
subgroup whose orbit is $\Sigma$ or term-by-term from (\ref{PQR}) now regarded
as a series of equations for the coefficients of this power series with $P,Q,R$
fixed). \end{proof}

The criteria in this theorem may be employed as follows. According to
\S\ref{normalforms} and especially Theorem~\ref{normforms}, the defining
equation of a homogeneous hypersurface with isotropy may be normalised to third
order. It is possible that all cubic terms vanish in which case we may consider
the consequences of Theorem~\ref{criteria} for $f(x,y,z)=2xy+z^2$ with~$N=3$.
Otherwise, we can take $f(x,y,z)$ to be one of
$$\begin{array}{ll}\mbox{\bf \rom{I3}}\quad 2xy+z^2+x^3&
                   \mbox{\bf \rom{I2}}\quad 2xy+z^2+x^2z\\[5pt]
                   \mbox{\bf \rom{I1}}\quad 2xy+z^2+x^2y-2xz^2\qquad&
                   \mbox{\bf \rom{I0}}\quad 2xy+z^2+3xyz-z^3.\end{array}$$
By way of illustration, let us consider in detail the case~{\bf I1} which is of
medium difficulty. There are several computations carried out with the aid of
{\sc Maple}. Further details on this use of computer algebra will be given
shortly.

The equations (\ref{PQR}) are polynomial in $x,y,z$ and so each coefficient
must vanish separately. In addition, since we are searching for a hypersurface
admitting
\begin{equation}\label{genI1}\mbox{\begin{mat}{cccc}
      0&0&0&0\\ 0&2&0&0\\ 0&0&1&0\\ 0&0&0&2\end{mat}}\end{equation}
as a generator of isotropy, we may normalise $P,Q,R$ by supposing that
\begin{equation}\label{normalisePQR}
p_{2,2}=0\quad q_{2,2}=0\quad r_{2,2}=0.\end{equation}
Altogether, this gives a system of linear equations for the entries
of $P,Q,R$ which is easily solved:
$$\begin{array}lP=\mbox{\begin{mat}{cccc}
2p_{3,3}-3&0&p_{1,3}&p_{1,4}\\
0&0&p_{2,3}&p_{2,4}\\
-p_{2,3}&-p_{1,3}&p_{3,3}&p_{3,4}\\
0&2&0&2p_{3,3}-2\end{mat}}\qquad
Q=\mbox{\begin{mat}{cccc}
2q_{3,3}&0&q_{1,3}&q_{1,4}\\
-1/2&0&q_{2,3}&q_{2,4}\\
-q_{2,3}&-q_{1,3}&q_{3,3}&q_{3,4}\\
2&0&0&2q_{3,3}\end{mat}}\\[30pt]
R=\mbox{\begin{mat}{cccc}
2r_{3,3}&0&r_{1,3}&r_{1,4}\\
0&0&r_{2,3}&r_{2,4}\\
2-r_{2,3}&-r_{1,3}&r_{3,3}&r_{3,4}\\
0&0&2&2r_{3,3}\end{mat}}\end{array}
$$
leaving the 18 entries
\begin{equation}\label{unknowns}
\begin{array}{llllll}p_{1,3},&p_{1,4},&q_{1,3},&q_{1,4},&r_{1,3},&r_{1,4},\\
  p_{2,3},&p_{2,4},&q_{2,3},&q_{2,4},&r_{2,3},&r_{2,4},\\
  p_{3,3},&p_{3,4},&q_{3,3},&q_{3,4},&r_{3,3},&r_{3,4}\end{array}\end{equation}
yet unknown. Now the first equation of (\ref{closure}) says that
$$\textstyle \frac{1}{2}(8p_{2,4}-5+4p_{3,3})x^2
+4(2p_{1,4}-q_{3,3}-2q_{2,4})xy-4q_{1,4}y^2+\cdots
+(7q_{1,3}+\cdots+2q_{2,3}r_{1,3})z^3$$
vanishes. Immediately, the coefficient of $y^2$ forces~$q_{1,4}=0$. More
specifically, if
\begin{equation}\label{startI1}f(x,y,z)=2xy+z^2+x^2y-2xz^2\end{equation}
can be completed to a power series
$F(x,y,z)$ defining a homogeneous hypersurface with
isotropy, then any normalised $Q$ will have $q_{1,4}=0$. This will eventually
be a consequence of the higher order terms and the
normalisation~(\ref{normalisePQR}). The coefficients of (\ref{closure}) give 41
polynomial constraints on~(\ref{unknowns}). There are just two solutions,
namely:
$$P=\mbox{\begin{mat}{cccc}  5/2    & 0   & 0     & 0      \\
                             0      & 0   & 0     & -3/4   \\
                             0      & 0   & 11/4  & 0      \\
                             0      & 2   & 0     & 7/2
          \end{mat}}\quad
Q=\mbox{\begin{mat}{cccc}    0      & 0   & 0     & 0      \\
                             -1/2   & 0   & 0     & 0      \\
                             0      & 0   & 0     & 0      \\
                             2      & 0   & 0     & 0
          \end{mat}}\quad
R=\mbox{\begin{mat}{cccc}    0      & 0   & 0     & 0      \\
                             0      & 0   & -1/2  & 0      \\
                             5/2    & 0   & 0     & 0      \\
                             0      & 0   & 2     & 0
          \end{mat}}$$
and
$$P=\mbox{\begin{mat}{cccc}  5/2    & 0   & 0     & 0      \\
                             0      & 0   & 0     & 1/2    \\
                             0      & 0   & 1/4   & 0      \\
                             0      & 2   & 0     & -3/2
          \end{mat}}\quad
Q=\mbox{\begin{mat}{cccc}    0      & 0   & 0     & 0      \\
                             -1/2   & 0   & 0     & 0      \\
                             0      & 0   & 0     & 0      \\
                             2      & 0   & 0     & 0
          \end{mat}}\quad
R=\mbox{\begin{mat}{cccc}    0      & 0   & 0     & 0      \\
                             0      & 0   & 2     & 0      \\
                             0      & 0   & 0     & 0      \\
                             0      & 0   & 2     & 0
          \end{mat}}$$
These give rise to cases~{\bf I1.1} and~{\bf I1.2} of
Theorem~\ref{maintheorem}. Specifically, if we add a general quartic term
$$f(x,y,z)=2xy+z^2+x^2y-2xz^2+\sum_{i+j+k=4}c_{i,j,k}x^iy^jz^k$$
and now re-consider (\ref{PQR}) with~$N=4$ and $P,Q,R$ one of these two
solutions of~(\ref{closure}), then the quartic terms are determined.
In fact, it is clear by inspection that (\ref{PQR}) with $N=4$ determines the
quartic terms which only enter the right hand sides. The crucial observation,
however, is that this overdetermined system is consistent as a consequence
of~(\ref{closure}). More precisely, the interpretation of (\ref{closure}) as
the closure of a subalgebra of the Lie algebra of affine motions implies that
the entire power series $F(x,y,z)$ may be defined implicitly by
\begin{equation}\label{exp}
\left\lgroup\exp\mbox{\begin{mat}{ccccc}
&&&&r\\
&\raisebox{-6pt}[0pt][0pt]{\makebox[0pt]{\hspace*{18pt}$rP+sQ+tR$}}&&&s\\
&&&&t\\  &&&&0\\
0\,&0\,&0\,&0\,&0\end{mat}}\right\rgroup
\mbox{\begin{mat}c0\\ 0\\ 0\\ 0\\ 1\end{mat}}=
\mbox{\begin{mat}cx\\ y\\ z\\ F(x,y,z)\\ 1\end{mat}}.\end{equation}
Complete details for the analogous case of surfaces are in~\cite[\S 2]{ee}.

To summarise then, cases~{\bf I1.1} and~{\bf I1.2} of Theorem~\ref{maintheorem}
are the only possible completions of (\ref{startI1}) defining a homogeneous
hypersurface with isotropy (necessarily generated by~(\ref{genI1})). It only
remains to check that these hypersurfaces really do have this isotropy. For
this, it suffices to take their $4^{\mathrm{th}}$ order truncations and apply
Theorem~\ref{criteria} with $N=4$ (without imposing the
normalisations~(\ref{normalisePQR})). It turns out that $P,Q,R$ are now
determined by (\ref{PQR}) alone up to adding arbitrary multiples
of~(\ref{genI1}). Furthermore, (\ref{closure}) now holds. So this constitutes
the full symmetry algebra and, apart from finding the explicit defining
functions given in Theorem~\ref{maintheorem}, cases~{\bf I1.1} and~{\bf I1.2}
are complete. Notice that, because Theorem~\ref{criteria} applies directly
when $N=4$, the higher order terms are uniquely determined simply by requiring
the hypersurface to be homogeneous irrespective of whether it has isotropy.
Finding explicit defining functions will be delayed until~\S\ref{compare}.

Though we can analyse all other cases in exactly the same way, there are some
initial observations which almost immediately deal with some of them. Take,
for example, the case~{\bf I2}. According to Theorem~\ref{normforms}, the only
possibility for isotropy in this case is the scaling~(\ref{scaling})
with~$t=2$. This limits the quartic terms to $bx^4$ for some $b$ whilst all
higher order terms must vanish. It is now easy to check that this equation does
indeed define a homogeneous hypersurface with this isotropy. Moreover, since
$b$ is unaffected by the only residual co\"ordinate freedom (namely, the
isotropy), it is a true parameter.

In case~{\bf I3},
it may be that (\ref{scaling}) with $t=3$ survives in the
isotropy of a corresponding hypersurface~$\Sigma$.
Straightaway this eliminates all
terms higher than cubic and we have case~{\bf I3} of Theorem~\ref{maintheorem}.
However, there remains the possibility that (\ref{scaling}) does not survive in
the isotropy of~$\Sigma$. This makes the isotropy one-dimensional, generated by
the third matrix of (\ref{generators}) with~$t=0$. The corresponding
one-parameter subgroup consists of null rotations~(\ref{nr}) and in
Theorem~\ref{maintheorem} we denote this case by~{\bf Inr}.

The detailed completion of cases~{\bf Inr} and~{\bf I0} follows the treatment
of {\bf I1} as above. The only real difficulty is in analysing the
criteria~(\ref{closure}) of Theorem~\ref{criteria}.
To ensure that all solutions of this system are found we
employed Buchberger's algorithm for Gr\"obner bases as implemented in the
`grobner' package of {\sc Maple} (Version V Release~3). In searching for
homogeneous hypersurfaces with scaling isotropy (\ref{scaling}) we can use
(\ref{normalisePQR}) but in case~{\bf Inr} we use $p_{2,3}=q_{2,3}=r_{2,3}=0$
instead. The entire analysis, including the calculation of
$4^{\mathrm{th}}$~order and (if necessary) $5^{\mathrm{th}}$~order terms, can
be completely automated. A {\sc Maple} program is available by
anonymous~ftp\footnote{%
ftp://ftp.maths.adelaide.edu.au/pure/meastwood/maple/{\bf thm2proof}}.
Unlike~{\bf I1}, in most other cases the closure
equations~(\ref{closure}) have infinitely many solutions with some entries in
$P,Q,R$ remaining free. These free entries show up in the higher order terms of
the corresponding completions as potential parameters. Having used the program
{\bf thm2proof} to find possible completions, there are two remaining tasks:
\begin{itemize}
\item apply the remaining co\"ordinate freedoms from Theorem~\ref{normforms}
to see whether they can be used to eliminate some of the parameters
appearing in these possible completions;
\item verify, by reapplying Theorem~\ref{criteria} with $N=4$ or~$5$, that
these completions really do have the anticipated isotropy.
\end{itemize}
There are two cases when the first of these tasks is non-trivial. When there
are no cubic terms, the typical output from {\bf thm2proof} is
$$\textstyle 2xy+z^2-2p_{1,4}x^2y^2-2p_{1,4}xyz^2-\frac{1}{2}p_{1,4}z^4.$$
(The output can vary depending on the particular invocation of {\sc Maple}
because the ordering it uses for computing Gr\"obner bases etcetera depends on
the internal addresses of the variables involved. This randomising effect can
be used to advantage by running the program several times and choosing the
simplest answer.) The rescaling (\ref{rescaling}) corresponding to the first
of~(\ref{generators}), has the effect of multiplying the coefficient of
$z^4$ by~$\lambda^2$. If this coefficient is non-zero, we may therefore
normalise it to unity and obtain case~{\bf Sp}. On the other hand, if it is
zero then we obtain case~{\bf Qd}. The other cases requiring special attention
are {\bf I3} and {\bf Inr} when $x^3$ is the only cubic term. No matter what
isotropy is assumed, the only possible quartic term is a multiple of~$x^4$.
A non-zero multiple may be normalised to $x^4$ itself by a suitable rescaling
(\ref{scaling}) with~$t=3$. This cuts down the residual co\"ordinate freedom to
(\ref{nr}) and leads to case~{\bf Inr}. When there are no quartic terms, then
we are led to~{\bf I3}. The final task of verifying that these hypersurfaces
really do have isotropy and, indeed, computing the full symmetry algebra and
checking that it closes is accomplished with a separate {\sc Maple}
program\footnote{%
ftp://ftp.maths.adelaide.edu.au/pure/meastwood/maple/{\bf thm2verify}}.

\section{Proof of Theorem~\protect\ref{explicit}}\label{compare}
There are two possible ways to proceed. We could start with the list of
explicit defining functions, verify that each of them gives a homogeneous
surface with isotropy, and then execute the normalisations of
\S\ref{normalforms} to obtain a perfect match with Theorem~\ref{maintheorem}.
We have written a {\sc Maple} program\footnote{%
ftp://ftp.maths.adelaide.edu.au/pure/meastwood/maple/{\bf thm1}} which takes a
defining function, computes its prospective symmetry algebra (by truncating its
power series as in Theorem~\ref{criteria}), checks that this algebra closes,
and then determines whether the hypersurface is genuinely invariant under these
symmetries (to infinite order). The 20 possibilities of the list are already in
{\bf thm1} and, indeed, this program shows them to be homogeneous with
isotropy. Of course, this approach is somewhat unsatisfactory because it does
not explain where the list comes from nor why minor variations such as
$$\begin{array}{lll}W=XY+Z^2\log Z\quad&WZ =XY+e^Z&WZ=XY+\log Z\\[5pt]
W^2=XY+\log Z\quad&W^2=XY+Z\log Z\quad&W^2=XY+Z^2\log Z\end{array}$$
are omitted. (According to Theorem~\ref{explicit}, they would already be on
the list but for a different choice of co\"ordinates.) In fact, none of these
is homogeneous as {\bf thm1} readily verifies. For example, the first of them
truncated at $4^{\mathrm{th}}$~order defines a closed symmetry algebra but does
not satisfy this algebra at $5^{\mathrm{th}}$~order. Rather,
$$\textstyle W+\frac{1}{28}+\frac{2}{7}Z=
               XY+\frac{9}{28}(\frac{5}{3}Z-\frac{2}{3})^{12/5}$$
is a homogeneous surface (\nn{7} with an affine change of co\"ordinates) which
just happens to have the same power series expansion about the point
$(0,0,0,1)$ up to $4^{\mathrm{th}}$~order.

More satisfactory is to start with Theorem~\ref{maintheorem} and derive
explicit defining functions in each case. The matrices $P,Q,R$ supplied by
{\bf thm2verify} describe the hypersurfaces parametrically (\ref{exp}) and,
with sufficient diligence, it is possible explicitly to solve for $F(x,y,z)$
and, after a suitable change of co\"ordinates, check the comparison table given
in~\S\ref{results}. There are, however, some observations which greatly
simplify this task. Cases~{\bf Qd} and~{\bf Sp} are clear by inspection (since
they are manifestly homogeneous and have the correct power series expansion up
to $4^{\mathrm{th}}$~order). Cases~{\bf I3} and~{\bf I2} are also immediate: it
was already observed in \S\ref{proofs} that their defining functions must be
polynomial.

In cases~{\bf I1} having (\ref{genI1}) generating the isotropy forces the
defining function $F(x,y,z)$ to have the form
$$F(x,y,z)=f(x)y+g(x)z^2.$$
The output from {\bf thm2verify} has
$$Q=\mbox{\begin{mat}{cccc}0&0&0&0\\
                           -1/2&0&0&0\\
                           0&0&0&0\\
                           2&0&0&0\end{mat}}\;+\mbox{ isotropy}$$
in both {\bf I1} cases. That the corresponding vector field
$$\left(1-\frac{x}{2}\right)\frac{\partial}{\partial y}
  +2x\frac{\partial}{\partial w}$$
be a symmetry implies that $f(x)=4x/(2-x)$. The vector fields corresponding to
$R$ distinguish {\bf I1.1} from~{\bf I1.2}:
$$2z\frac{\partial}{\partial y}
  +\frac{\partial}{\partial z}
  +2z\frac{\partial}{\partial w}\quad\mbox{ versus }\quad
  -\frac{z}{2}\frac{\partial}{\partial y}
  +\left(1+\frac{5}{2}x\right)\frac{\partial}{\partial z}
  +2z\frac{\partial}{\partial w}.$$
They determine $g(x)$ as $(2-5x)/(2-x)$ or $4/((2-x)(2+5x))$, respectively.
These are the defining functions given in Theorem~\ref{maintheorem} and
the affine co\"ordinate changes
$$\mbox{\begin{mat}cx\\ y\\ z\\ w\end{mat}}
 =\mbox{\begin{mat}{cccc}0 &-2/5&0 &0\\
                         -1&0   &-5&0\\
                         0 &0   &0 &2\\
                         4 &0   &0 &0\end{mat}}
\mbox{\begin{mat}cW\\ X\\ Y\\ Z\end{mat}}
$$and$$
\mbox{\begin{mat}cx\\ y\\ z\\ w\end{mat}}
 =\mbox{\begin{mat}{cccc}0  &2/5&0 &0 \\
                         2  &-2 &-7&-6\\
                         2  &-2 &0 &0 \\
                         -8 &8  &8 &4\end{mat}}
\mbox{\begin{mat}cW\\ X\\ Y\\ Z\end{mat}}
+\mbox{\begin{mat}c-2/5\\ 6\\ 0\\ -4\end{mat}}$$
give the defining functions~\nn{5} and~\nn{6} with their respective basepoints.

Case~{\bf Inr} has
$$R=\mbox{\begin{mat}{cccc}0&0&0&0\\ 0&0&0&0\\ 0&0&0&0\\ 0&0&2&0\end{mat}}
\;+\mbox{ isotropy}$$
which implies that
$$F(x,y,z)=z^2+\theta(x,y)$$
and it follows easily that $\{w=\theta(x,y)\}$ is a homogeneous surface in
affine three-space. These were classified in \cite{dkr,ee,lob} and in \cite{ee}
a method was given for locating any given surface. In fact, it is easy to spot
that these surfaces are exactly the class~{\bf N6} of \cite[Theorem~2]{ee} with
almost the same normalisation: the parameter~$b$ is exactly as in~\cite{ee}.
Furthermore, in \cite[\S6.2]{ee} was given a precise comparison between these
normal forms and the explicit defining functions of~\cite{dkr}. Following this
through gives \nn{16}--\nn{20} as in the comparison table of~\S\ref{results}.
(Though in \cite{dkr} the link between symmetry algebra~\#9 and surface~\#12
should have $\alpha$ replaced by $1/\alpha$.)

There is a similar link with homogeneous surfaces in cases~{\bf I0}. The
isotropy implies that
$$F(x,y,z)=\theta(u,z)\quad\mbox{where }u=xy$$
whence
$$\begin{array}l
\begin{array}l\mbox{\begin{mat}c\displaystyle \frac{\partial F}{\partial x},
  \frac{\partial F}{\partial y},\frac{\partial F}{\partial z},-1\end{mat}}
\left[\mbox{\begin{mat}{cccc}p&0&0&0\\
                             0&q&0&0\\
                             0&0&a&b\\
                             0&0&c&d\end{mat}}
\mbox{\begin{mat}cx\\ y\\ z\\ F\end{mat}}+
\mbox{\begin{mat}c0\\ 0\\ r\\ s\end{mat}}\right]\\[20pt]
\hspace*{140pt}=
\mbox{\begin{mat}c\displaystyle \frac{\partial \theta}{\partial u},
                                \frac{\partial \theta}{\partial z},-1\end{mat}}
\left[\mbox{\begin{mat}{ccc}p+q&0&0\\ 0&a&b\\ 0&c&d\end{mat}}
\mbox{\begin{mat}cu\\ z\\ \theta\end{mat}}+
\mbox{\begin{mat}c0\\ r\\ s\end{mat}}\right]\end{array}\\[45pt]
\begin{array}l\mbox{\begin{mat}c\displaystyle \frac{\partial F}{\partial x},
  \frac{\partial F}{\partial y},\frac{\partial F}{\partial z},-1\end{mat}}
\left[\mbox{\begin{mat}{cccc}0&0&r&s\\
                             0&0&0&0\\
                             0&p&0&0\\
                             0&q&0&0\end{mat}}
\mbox{\begin{mat}cx\\ y\\ z\\ F\end{mat}}+
\mbox{\begin{mat}ct\\ 0\\ 0\\ 0\end{mat}}\right]\\[20pt]
\hspace*{140pt}=y
\mbox{\begin{mat}c\displaystyle \frac{\partial \theta}{\partial u},
                                \frac{\partial \theta}{\partial z},-1\end{mat}}
\left[\mbox{\begin{mat}{ccc}0&r&s\\ 0&0&0\\ 0&0&0\end{mat}}
\mbox{\begin{mat}cu\\ z\\ \theta\end{mat}}+
\mbox{\begin{mat}ct\\ p\\ q\end{mat}}\right]\end{array}\\[45pt]
\begin{array}l\mbox{\begin{mat}c\displaystyle \frac{\partial F}{\partial x},
  \frac{\partial F}{\partial y},\frac{\partial F}{\partial z},-1\end{mat}}
\left[\mbox{\begin{mat}{cccc}0&0&0&0\\
                             0&0&r&s\\
                             p&0&0&0\\
                             q&0&0&0\end{mat}}
\mbox{\begin{mat}cx\\ y\\ z\\ F\end{mat}}+
\mbox{\begin{mat}c0\\ t\\ 0\\ 0\end{mat}}\right]\\[20pt]
\hspace*{140pt}=x
\mbox{\begin{mat}c\displaystyle \frac{\partial \theta}{\partial u},
                                \frac{\partial \theta}{\partial z},-1\end{mat}}
\left[\mbox{\begin{mat}{ccc}0&r&s\\ 0&0&0\\ 0&0&0\end{mat}}
\mbox{\begin{mat}cu\\ z\\ \theta\end{mat}}+
\mbox{\begin{mat}ct\\ p\\ q\end{mat}}\right].\end{array}
\end{array}$$
It turns out that the matrices $P,Q,R$ in cases~{\bf I0} have the forms
indicated in these identities. For example, in case~{\bf I0.1} with $b=6$
{\bf thm2verify} gives
$$P=\mbox{\begin{mat}{cccc}0&0&0&0\\
               0&0&0&0\\ 0&-3/2&0&0\\ 0&2&0&0\end{mat}}
\quad
Q=\mbox{\begin{mat}{cccc}0&0&0&0\\
                0&0&0&0\\ -3/2&0&0&0\\ 2&0&0&0\end{mat}}
\quad
R=\mbox{\begin{mat}{cccc}15/2&0&0&0\\
             0&0&0&0\\ 0&0&6&-9/8\\ 0&0&2&9\end{mat}}$$
apart from isotropy. Therefore, the surface $\{w=\theta(u,z)\}$ is homogeneous
and, following the notation of~\cite{dkr}, its symmetry algebra contains
\begin{equation}\label{syms}
\mbox{\begin{mat}{ccc}15/2&0&0\\ 0&6&-9/8\\ 0&2&9\end{mat}}
  +\mbox{\begin{mat}c0\\ 0\\ 1\end{mat}}
\quad\mbox{ and }\quad
  \mbox{\begin{mat}{ccc}0&0&0\\ 0&0&0\\ 0&0&0\end{mat}}
  +\mbox{\begin{mat}c1\\ -3/2\\ 2\end{mat}}.\end{equation}
The second of these generates a uniform translation so the surface is a
cylinder, i.e.\ class {\bf D2} of~\cite{ee}. {From} the full symmetry algebra
of {\bf D2} given in \cite{ee} it is easy to check that the true
parameter $a^2$, if non-zero, is given by
$$\frac{32[3\trace(M)\trace(M^2)-2\trace(M^3)
     -3\trace(M^2)\lambda-3\trace(M)\lambda^2+5\lambda^3]}
       {25[\trace(M)-\lambda]
       [5\trace(M)^2-9\trace(M^2)-10\trace(M)\lambda+14\lambda^2]}$$
for any non-zero $M$ from the matrix part of the algebra where $\lambda$ is the
eigenvalue of $M$ for the translation vector.
Thus $a^2=64/25$ and from \cite[\S1]{ee} the surface must be
$\{Z=X\log X\}$ for a suitable choice of affine co\"ordinates.

Unfortunately, this abstract reasoning loses track of the distinguished
co\"ordinate $u=xy$ so a more direct argument must be employed. The change of
co\"ordinates
$$\mbox{\begin{mat}cu\\ z\\ w\end{mat}}\longmapsto
\mbox{\begin{mat}{ccc}1&0&0\\ 0&-3/2&1\\ 0&2&0\end{mat}}
\mbox{\begin{mat}cu\\ z\\ w\end{mat}}$$
preserves $u$ but conjugates the symmetries (\ref{syms}) to
$$\mbox{\begin{mat}{ccc}15/2&0&0\\ 0&15/2&1\\ 0&0&15/2\end{mat}}
  +\mbox{\begin{mat}c0\\ 1/2\\ 3/4\end{mat}}
\quad\mbox{ and }\quad
  \mbox{\begin{mat}{ccc}0&0&0\\ 0&0&0\\ 0&0&0\end{mat}}
  +\mbox{\begin{mat}c1\\ 1\\ 0\end{mat}}.$$
Now we can employ the surface version of~(\ref{exp}):
$$\begin{array}l\left\lgroup\exp\mbox{\begin{mat}{cccc}
15s/2 & 0     & 0 & t \\
0     & 15s/2 & s & s/2 + t\\
0     & 0     & 15s/2 &3s/4\\
0 & 0 & 0 & 0\end{mat}}\right\rgroup
\mbox{\begin{mat}c 0\\ 0\\ 0\\ 1\end{mat}}\\[30pt]
\hspace*{150pt}= \mbox{\begin{mat}c\frac{2}{15}t(e^{15s/2}-1)/s\\
\frac{1}{10}se^{15s/2}+\frac{2}{75}(2s+5t)(e^{15s/2}-1)/s \\
\frac{1}{10}(e^{15s/2}-1)\quad \\ 1\end{mat}}
=\mbox{\begin{mat}cu\\ z\\ w\\ 1\end{mat}}.\end{array}$$
This may be solved:
$$75z=75u+(1+10w)\log(1+10w)+40w$$
and a further affine change of co\"ordinates
$$W=75z-40w \qquad X=75x \qquad Y=y \qquad Z=1+10w$$
evidently gives \nn{10}. All other cases follow similarly.

\section{Remarks on the real case}\label{remarks}
The analysis in the real case proceeds as for the complex case save for a few
minor changes. Notice that in \S\ref{proofs}, only twice was it used that we
were working over the complex numbers. It was when we were normalising the
quartic terms in cases {\bf Sp} and~{\bf Inr}. Taking this into account, it
follows that these cases have two real forms:
\begin{description}\addtolength{\itemsep}{-7pt}
                   \setlength{\labelwidth}{30pt}
                   \setlength{\itemindent}{0pt}
\item[\rom{Sp}${}^{\pm}$\hfill]
                        $w=2xy+z^2\pm 4x^2y^2\pm 4xyz^2\pm z^4+{\mathrm{O}}(5)
                          =\pm(1-\sqrt{1\mp 4(2xy+z^2)})/2$
\item[\rom{Inr}${}^{\pm}$\hfill]  $w=2xy+z^2+x^3\pm x^4+bx^5+{\mathrm{O}}(6)$
\end{description}
Therefore, the conclusion to be drawn from \S\ref{proofs} in the real case is
that, with these two exceptions, if we start off a real power series as
$w=2xy+z^2+\cdots$ and continue with no cubic terms or with cubic terms listed
in Theorem~\ref{normforms}, then the resulting list of real homogeneous
hypersurfaces with isotropy is just as in Theorem~\ref{maintheorem}.

However, before arriving at \S\ref{proofs} we were normalising the defining
equations in~\S\ref{normalforms} and here also, complex numbers entered at two
stages. The first was in normalising the quadratic terms. Over the reals there
are two possibilities, namely
$$w=2xy+z^2+{\mathrm O}(3)\quad\mbox{ and }
  \quad w=x^2+y^2+z^2+{\mathrm O}(3),$$
{\em hyperbolic\/} and {\em elliptic}. This now shows up in the second stage
where we conjugated and rescaled a matrix in ${\frak o}(3,{\Bbb C})$ into one
of the three standard forms~(\ref{threestandardforms}). For ${\frak o}(2,1)$
there are four standard forms because non-null vectors now come in two
flavours, either space-like or time-like. Over ${\mathrm O}(3,{\Bbb C})$ the
corresponding matrices are conjugate up to scale:
$$\sqrt{2}i\mbox{\begin{mat}{ccc}i/\sqrt{2}&i/\sqrt{2}&1\\
                      i/\sqrt{2}&i/\sqrt{2}&-1\\
                      1&-1&0\end{mat}}
  \mbox{\begin{mat}{ccc}-1&0&0\\ 0&1&0\\ 0&0&0\end{mat}}
\mbox{\begin{mat}{ccc}i/\sqrt{2}&i/\sqrt{2}&1\\
                      i/\sqrt{2}&i/\sqrt{2}&-1\\
                      1&-1&0\end{mat}}^{-1}
  =\mbox{\begin{mat}{ccc}0&0&1\\ 0&0&1\\ -1&-1&0\end{mat}}$$
but not over~${\mathrm O}(2,1)$. As a corresponding complex co\"ordinate
change we may choose
$$\mbox{\begin{mat}cx\\ y\\ z\end{mat}}\longmapsto
\mbox{\begin{mat}{ccc}i/\sqrt{2}&i/\sqrt{2}&1\\
                      i/\sqrt{2}&i/\sqrt{2}&-1\\
                      1&-1&0\end{mat}}
\mbox{\begin{mat}cx\\ y\\ z\end{mat}}\qquad
w\mapsto-2w$$
which takes $w=2xy+z^2+3xyz-z^3+\cdots$ to
$$\textstyle w=2xy+z^2
  +\frac{5}{4}x^3-\frac{3}{4}(x^2y-2xz^2)
  +\frac{3}{4}(xy^2-2yz^2)-\frac{5}{4}y^3+\cdots$$
and gives two real forms for each of~{\bf I0}. This is the full extent of this
alternative normalisation. So cases~{\bf Sp}, {\bf Inr}, and {\bf I0} in
Theorem~\ref{maintheorem} have two hyperbolic real forms and the rest have just
one.

It is easy to check that the Pick invariant is zero in cases~{\bf I3},
{\bf I2}, and~{\bf I1} of Theorem~\ref{normforms}. Therefore, none of the
corresponding hypersurfaces in Theorem~\ref{maintheorem} (including case~{\bf
Inr}) can have an elliptic real form. The complex change of co\"ordinates
$$\mbox{\begin{mat}cx\\ y\end{mat}}\longmapsto
\mbox{\begin{mat}{cc}1/\sqrt{2}&i/\sqrt{2}\\ 1/\sqrt{2}&-i/\sqrt{2}\end{mat}}
\mbox{\begin{mat}cx\\ y\end{mat}}$$
gives the unique elliptic form of case~{\bf Qd} and cases~{\bf I0}. It gives
one of the two real forms of~{\bf Sp}, namely the sphere. The other is the
hyperhyperboloid of two sheets. Assembling these observations and tracing
through to the explicit defining functions gives the following real
classification list.
\renewcommand{\arraystretch}{1.2}
\begin{center}\begin{tabular}{|l|ll|}                         \hline
\nn{1} &$W=XY+Z^2$                &$W=X^2+Y^2+Z^2$          \\ \hline
\nn{2} &$W^2=XY\pm Z^2+1$         &$\pm W^2=X^2+Y^2+Z^2\pm 1$\\ \hline
\nn{3} &$W=XY+Z^2+X^3$            &                         \\ \hline
\nn{4} &$W=XY+Z^2+X^2Z+\alpha X^4$&                         \\ \hline
\nn{5} &$W=XY+Z^2+XZ^2$           &                         \\ \hline
\nn{6} &$W^2=XY+X^2Y+X^2Z$        &                         \\ \hline
\nn{7} &$W=XY+Z^\alpha$           &$W=X^2+Y^2\pm Z^\alpha$  \\ \hline
\nn{8} &$W=XY+e^Z$                &$W=X^2+Y^2\pm e^Z$       \\ \hline
\nn{9} &$W=XY+\log Z$             &$W=X^2+Y^2\pm \log Z$    \\ \hline
\nn{10}&$W=XY+Z\log Z$            &$W=X^2+Y^2\pm Z\log Z$   \\ \hline
\nn{11}&$W^2=XY+Z^\alpha$         &$W^2=X^2+Y^2\pm Z^\alpha$\\ \hline
\nn{12}&$W^2=XY+e^Z$              &$W^2=X^2+Y^2\pm e^Z$     \\ \hline
\nn{13}&$WZ=XY+Z^\alpha$          &$WZ=X^2+Y^2\pm Z^\alpha$ \\ \hline
\nn{14}&$WZ=XY+Z\log Z$           &$WZ=X^2+Y^2\pm Z\log Z$  \\ \hline
\nn{15}&$WZ=XY+Z^2\log Z$         &$WZ=X^2+Y^2\pm Z^2\log Z$\\ \hline
\nn{16}&$W=XY\pm Z^2+X^\alpha$    &                         \\ \hline
\nn{17}&$W=XY\pm Z^2+e^X$         &                         \\ \hline
\nn{18}&$W=XY\pm Z^2+\log X$      &                         \\ \hline
\nn{19}&$W=XY\pm Z^2+X\log X$     &                         \\ \hline
\nn{20}&$W=XY\pm Z^2+X^2\log X$   &                         \\ \hline
\end{tabular}\end{center}\renewcommand{\arraystretch}{1}

\small\renewcommand{\section}{\subsection}

{\sc \begin{flushleft}
Department of Pure Mathematics\\
University of Adelaide\\
South Australia 5005\\[5pt]
\normalshape
E-mail: meastwoo@maths.adelaide.edu.au\\
\mbox{\phantom{E-mail: }vezhov@maths.adelaide.edu.au}
\end{flushleft}}
\end{document}